# HIGHER TOPOLOGICAL COMPLEXITY OF ARTIN TYPE GROUPS

SERGEY YUZVINSKY

ABSTRACT. We calculate the higher topological complexity $TC_s$ for the complements of reflection arrangements, in other words for the pure Artin type groups of all finite complex reflection groups. In order to do that we introduce a simple combinatorial criterion of arrangements sufficed for the cohomological low bound for $TC_s$ to coincide with the dimensional upper bound.

## 1. INTRODUCTION

Topological complexity of a topological space $X$ ($TC(X)$) was defined by M. Farber in [6] as a specialization of the Schwarz genus [12]. Unlike the Schwarz genus in general, $TC(X)$ is an invariant of the homotopy type of $X$. Later Yu. Rudyak in [11] extended Farber's definition to higher topological complexity $TC_s(X)$ for $s = 2, 3, 4, \ldots$ such that $TC(X)=TC_2(X)$. One of common features of these invariants is a lower bound determined by the ring structure of $H^*(X)$.

This paper is concerned with the special class of topological spaces - the complements of a complex hyperplane arrangement. Previously $TC_2$ has been calculated for particular classes of arrangements such as Coxeter series in [7] and general position arrangements in [13, 4]. These examples prompted the Conjecture that for all arrangement complements $TC_2$ coincides with the cohomological low bound. The only known results for arrangement complements and arbitrary $s$ (besides the basic examples of the circle and tori in [11, 2]) is the calculation for the Coxeter series of type $A$ in the recent preprint [8].

In the present paper we give a simple combinatorial condition sufficed for the cohomological low bound to coincide with the dimensional upper bound. This allows us to compute $TC_s$ for a wide class of arrangements (including all complex reflection arrangements, i.e., $K[\pi, 1]$'s for the pure Artin type groups). In all arrangements of this class the value of $TC_s$ coincides with the cohomological low bound for this $s$.

The results of this paper constituted a talk given by the author at the conference "Configuration Spaces" in Cortona in September of 2014. The author is grateful to the organizers of the very useful and pleasant conference for inviting him.

## 2. DEFINITION OF $TC_s$ AND MAIN PROPERTIES

**Definition 2.1.** *Let $X$ be a path-connected topological space and $s$ an integer at least 2. Then $TC_s(X)$ is the Schwarz genus of the fibration*

$$\phi_s : X^{[0,1]} \to X^s,$$

*where*

$$\phi_s(\gamma) = \left(\gamma(0), \gamma\left(\frac{1}{s-1}\right), \gamma\left(\frac{2}{s-1}\right), \ldots, \gamma\left(\frac{s-2}{s-1}\right), \gamma(1)\right).$$

In other words, it is the smallest number $n$ such that $X^s$ is partitioned into Euclidean neighborhood retracts $W_i$ ($i = 0, 1, \ldots, n$) and on each $W_i$ there exists a section $f_i : X^s \to X^{[0,1]}$ of $\phi_s$ (i.e., $\phi_s \circ f_i = 1_{W_i}$). These data is called a *motion planning* (m.p.). We use the *reduced* (or normalized) version of TC such that $\mathrm{TC}_s(X) = 0$ for a contractible $X$ and each $s$.

Note that $\phi_s$ is a fibrational substitute with the fiber $(\Omega X)^{s-1}$ for the diagonal imbedding $d_s : X \to X^s$.

The following properties can be found in [11, 2].

(1) $\mathrm{TC}(X)$ is an invariant of the homotopy type of $X$.

(2) $\mathrm{TC}_s(X) \leq s \cdot \mathrm{hdim}(X)$ where hdim is the homotopy dimension (*the dimensional upper bound*).

(3) $\mathrm{TC}_s(X \times Y) \leq \mathrm{TC}_s(X) + \mathrm{TC}_s(Y)$ (*the product formula*).

(4) The cohomological lower bound.

This is the only low bound and it requires some definitions.

**Definition 2.2.** *Let $d_s$ be the diagonal embedding $X \to X^s$. Denote by $cl(X, s)$ the cup length in $\ker d_s^*$, i.e., the largest integer $k$ for which there exist $k$ elements $u_i \in H^*(X^s)$ such that $d_s^* u_i = 0$ for every $i$ and $u_1 u_2 \cdots u_k \neq 0$.*

We have the following *cohomological low bound*:

$$\mathrm{TC}_s(X) \geq cl(X, s).$$

This inequality holds for cohomology with arbitrary coefficients, even for local coefficients. In the rest of the paper we will use cohomology with coefficients in $\mathbb{C}$ omitting coefficients from the notation.

**Example**

$\mathrm{TC}_s(S^1) = s - 1$ for every $s$.

Indeed choose an orientation and denote by $u$ the generator of $H^1(S^1)$. Then the elements

$$u^{(i)} = u \otimes 1 \otimes \cdots \otimes 1 - 1 \otimes \cdots \otimes 1 \otimes u \otimes 1 \otimes \cdots \otimes 1$$

where the second $u$ is in the $i$-th position are in $\ker d^*$ and $\prod_{i=2}^{s} u^{(i)} \neq 0$. By the cohomological lower bound $\mathrm{TC}_s(S^1) \geq s - 1$.

For a m.p. one can use the covering of the torus $(S^1)^s$ by $D_k = \{(x_1, \ldots, x_s)\}$ ($k = 0, 1, \ldots, s-1$) such that $x_j = x_{j+1}$ for precisely $k$ indexes $j$. A path $\gamma_j(x)$ from $x_j$ to $x_{j+1}$ is constant if $x_j = x_{j+1}$ and the rotation along the fixed orientation of $S^1$ otherwise.

## 3. Complement of hyperplane arrangement

In this paper we will deal mostly with the topological spaces that are hyperplane arrangement complements..

**Definition 3.1.** *A (complex linear essential) hyperplane arrangement is a set $\mathcal{A}$ of $n$ linear hyperplanes in $\mathbb{C}^r$ such that $\bigcap_{H \in \mathcal{A}} H = \{0\}$. The arrangement complement is the topological space $M = \mathbb{C}^r \setminus \bigcup_{H \in \mathcal{A}} H$.*

Among the arrangement complements there are, for instance, $K[\pi, 1]$ where $\pi$ is the pure Artin type group for an arbitrary finite complex reflection group. The most frequently used examples of that are the Braid arrangements.

**Example.** Consider $n = \binom{\ell}{2}$ hyperplanes given by the equations $x_i = x_j$ for all $1 \leq i < j \leq \ell$. This arrangement is called Braid arrangement because $M$ is $K[\pi, 1]$ where $\pi = \pi_1(M)$ is the pure Braid group on $\ell$ strings, that is the pure Artin group for $\Sigma_\ell$.

For an arbitrary arrangement the algebra $A = H^*(M)$ is well-known from work of Arnold, Brieskorn and Orlik-Solomon ([1, 3, 9]).

For each hyperplane $H \in \mathcal{A}$ fix a linear form $\alpha_H$ with $\ker \alpha_H = H$. Then $A$ can be identified with the subalgbera of the algebra of all the holomorphic differential forms on $M$ generated by the logarithmic forms $\frac{d\alpha_H}{\alpha_H}$ ($H \in \mathcal{A}$). The classes $e_H$ of these forms form a canonical basis of $A^1$ whence for every $x \in A^1$ we have $x = \sum_{H \in \mathcal{A}} x_H e_H$ for some $x_H \in \mathbb{C}$. Relations for the generators are explicitly described and can be found in [10]. These relations imply in particular that $H^p(M) = 0$ for $p > r$.

A stronger fact is that $M$ has the homotopy type of a finite simplicial complex of dimension $r$ (see [10]).

The relations imply also that $A$ is determined by the combinatorics of $\mathcal{A}$, i.e., the collection of linearly independent subsets of $\mathcal{A}$ (called *simple matroid*). In particular the (square-free) monomials corresponding to dependent sets of hyperplanes vanish in $A$. Hence the (square-free) monomials corresponding to independent sets ("independent monomials") linearly generate $A$ but they are not linearly independent (over $\mathbb{C}$) in general. Theory of Gröbner basis gives so called no-broken-circuit (**nbc**) monomials that do form a basis of $A$.

To define this basis we need to fix a linear order on $\mathcal{A}$ whence on $\{e_1, \ldots, e_n\}$ which gives the deg-lex order on the monomials. Then a *circuit* is a minimal dependent set of $e_i$ and a *broken circuit* is circuit with the smallest element (in the fixed order) omitted. Then an **nbc**-monomial is a monomial whose support does not contain any broken circuits. It is easy to see that the set of **nbc**-monomials form the basis of $A$ given by the Gröbner theory for the deglex monomial ordering.

Later in this paper we will use the following.

**Property (\*) of nbc basis.**

Suppose an order is fixed on $\mathcal{A}$ and $\mu$ is a non-**nbc** monomial for this order. Then its representation as a linear combination of **nbc** monomials looks like

$$\mu = \sum_i \pm \mu_i$$

where for each **nbc**-monomial $\mu_i$ we have $\mu_i < \mu$ in the deglex order.

## 4. Properties of $\mathrm{TC}_s(M)$

1. The general upper bound for $M$ can be made a little more tight. Namely

$$\mathrm{TC}_s(M) \leq sr - 1.$$

Indeed for a non-empty arrangement $M = \bar{M} \times \mathbb{C}^*$ where $\bar{M}$ is the projectivization of $M$ that has the homotopy type of a CW-complex of dimension $r-1$. Hence by the product formula

$$\mathrm{TC}_s(M) \leq \mathrm{TC}_s(M_0) + \mathrm{TC}_s(S^1) \leq s(r-1) + s - 1 = sr - 1.$$

2. To find a lower bound we need some preparation.

### 4.1. Products over subsets.
Fix an integer $s \geq 2$ and for each generator $e_i \in H^1(M)$, and each $j$ $(1 < j \leq s)$ put

$$e_i^{(j)} = e_i \otimes 1 \otimes \cdots \otimes 1 - 1 \otimes \cdots \otimes 1 \otimes e_i \otimes 1 \otimes \cdots \otimes 1$$

where $e_i$ in the second summand is in the $j$th position. Clearly each $e_i^{(j)} \in \ker d_s^*$.

Then for every $I \subset \{e_i^{(j)}\}$ if the product of all elements from $I$ does not vanish then $|I|$ is a lower bound for $\mathrm{TC}_s$.

### 4.2. Products over pairs.
In the rest of the paper we will identify subsets of $\bar{n}$ with the respective subsets of generators in $A^1$ and with subarrangements of $\mathcal{A}$. The rank $\operatorname{rk} S$ of a subset $S$ is the rank of the respective subarrangement which is the cardinality of its base (i.e., a maximal independent set). The rank of $\bar{n}$ is $r$.

Let $Q = (B, C)$ be an ordered pair of disjoint subsets of $\bar{n}$. The *product over $Q$* is defined by the formula
$$\pi_Q = \pi_B \cdot \pi'_C$$
where
$$\pi_B = \prod_{i \in B} \prod_{j=2}^{s} e_i^{(j)}, \quad \pi'_C = \prod_{i \in C} e_i^{(2)}.$$

We put $\bar{Q} = B \cup C$.

### 4.3. Basic pairs and balanced sets.

**Definition 4.1.** *A pair $Q$ is basic if $|B| = r$ (i.e., $B$ is a base) and $B, C$ are **nbc** in $\bar{Q}$ for some linear order on it.*

**Remark.**
We can extend a linear order on $\bar{Q}$ to a linear order on $\mathcal{A}$ so that every element of $\bar{Q}$ are smaller than every element of $\mathcal{A} \setminus \bar{Q}$. Then any monomial with support in $\bar{Q}$ is **nbc** in $\bar{Q}$ if and only if it is **nbc** in the whole $\mathcal{A}$.

**Definition 4.2.** *A subset $S \subset \bar{n}$ of full rank is balanced if for any its (linearly) closed non-empty subsets $S'$ we have $|S'| < 2\operatorname{rk}(S')$.*

**Theorem 4.1.** *(i) If a pair $Q$ is basic then $\bar{Q}$ is balanced.*
*(ii) Every balanced set $S$ is $\bar{Q}$ for some basic pair $Q$, i.e., is the union of the elements of the pair.*

*Proof.* (i) Suppose a pair $Q = (B, C)$ is basic and fixed an order such that $B$ and $C$ are **nbc**. Also assume there exists a closed non-empty subset $D \subset \bar{Q}$ with $|D| \geq 2\operatorname{rk}(D)$. Since the sets $B \cap D$ and $C \cap D$ are independent they both are bases of $D$. Now if $i$ is the least element of $D$ then to whichever of two bases it belongs, it depends on the other base which contradicts to $B$ and $C$ being **nbc** in $\bar{Q}$.

(ii) Suppose a set $S$ is balanced and choose a base $B$ of $S$ hence of $\mathcal{A}$. Put $C = S \setminus B$. Since $S$ is balanced $s = |C| \leq r - 1$. Using that $S$ is balanced again we can find $r - s + 1$ elements in $B$ independent of $C$. Order them linearly from 1 to $r - s + 1$ and call the set they form by $B_1$. Again by the same property there exist two elements in $C$ independent of $B \setminus B_1$. Assign numbers $r - s + 2$ and $r - s + 3$ to these elements and call the set of them by $C_1$. Notice that $|C \setminus C_1| = s - 2$ and $|B \setminus B_1| = s - 1$. Now we just repeat the reasoning. There exist two elements in $B \setminus B_1$ independent on $C_1$ and we can assign numbers $r - s + 4$ and $r - s + 5$ to them. Continuing this process we obtain at some step a linear ordering on $S$ such that no element depends on the set of larger (in this ordering) elements. Thus $B$ and $C$ are **nbc** in $S$ for this order whence the pair is basic. $\square$

## 5. Calculation of a lower bound

**Theorem 5.1.** *Let $\mathcal{A}$ be a central arrangement. Then for every integer $s$, $s \geq 2$, and every basic pair $Q = (B, C)$ we have $\pi_Q \neq 0$. Hence $TC_s \geq (s-1)r + |C|$.*

*Proof.* By construction, $\pi_Q$ is the sum of pure tensors with coefficients $\pm 1$ among which there is $\mu = e_C \otimes e_B \cdots \otimes e_B$ where for every subset $S \subset [n]$ we put $e_S = \prod_{i \in S} e_i$. Since $Q$ is basic all monomials in $\mu$ are **nbc** in some order on $\mathcal{A}$ that we fix. Thus it suffices to proof that no other simple tensor from $\pi_Q$ contains $\mu$ in the decomposition of its monomials in the linear combinations of **nbc** monomials.

Suppose that $\nu = \nu_1 \otimes \nu_2 \otimes \cdots \otimes \nu_s$ is such a simple tensor. Since the monomials $\nu_j$ for $j > 2$ cannot have elements from $C$ and have degree $r$ then $\nu_i = e_B$ for $i > 2$. The first two monomials are $e_{C_i}$ ($i = 1, 2$) where $(C_1, C_2)$ is a partition of $B \cup C$ (with $|C_2| = r$).

Using Property (*) of the **nbc** basis we obtain the following. If $e_{C_1} \otimes e_{C_2}$ contains $e_C \otimes e_B$ in the decomposition and at least one of $e_{C_i}$ is not **nbc** then $e_{C_1} \geq e_C$ and $e_{C_2} \geq e_B$ with at least one of the inequalities is strict. This contradicts the fact that $(C_1, C_2)$ and $(C, B)$ are partitions of $C \cup B$. Thus $\nu = \mu$ whence $\mu$ cannot be cancelled. $\square$

## 6. Large arrangements

**Definition 6.1.** *We call an arrangement large if there exists a basic pair $(B, C)$ with $|C| = r - 1$.*

Comparing this with the dimensional upper bound for $M$ we obtain for large arrangements that
$$\mathrm{TC}_s(M) = sr - 1$$
for every $s$.

Large arrangements are easy to find due to the following sufficient condition.

**Definition 6.2.** *A pair $(B, C)$ is well-balanced if $B$ is a base, $|C| = r - 1$, and no $b \in B$ depends on $C$. An arrangement is well-balanced if there is a well-balanced pair in it.*

**Theorem 6.1.** *Every well-balanced pair is balanced.*

*Proof.* Indeed suppose $(B, C)$ is well-balanced but there is a non-empty $D \subset B \cup C$ with $|D| \geq 2\,\mathrm{rk}\,D$. Then $D \cap B$ and $D \cap C$ are independent whence both are bases of $D$. Hence every $b \in D \cap B$ depends of $D \cap C$ which contradicts the condition. $\square$

Let $L(\mathcal{A})$ be the lattice of all intersections of hyperplanes from $\mathcal{A}$ ordered opposite to inclusion. For $X \in L(\mathcal{A})$ we put $\mathcal{A}_X = \{H \in \mathcal{A} | H \geq X\}$.

**Definition 6.3.** $L(\mathcal{A})$ *is well-balanced if there exists* $X \in L(\mathcal{A})$, $\operatorname{codim} X = r - 1$ *such that for no* $Y \in L(\mathcal{A}) \setminus \{0\}$ *we have* $\mathcal{A} = \mathcal{A}_X \cup \mathcal{A}_Y$.

This definition makes sense for an arbitrary finite geometric lattice.

**Theorem 6.2.** *If* $L(\mathcal{A})$ *is well-balanced then there exists a well-balanced pair in* $\mathcal{A}$.

*Proof.* Let $C$ be a base of $\mathcal{A}_X$ from the definition. Then $|C| = r - 1$. Put $\mathcal{A}' = \mathcal{A} \setminus \mathcal{A}_X$. By definition $\operatorname{rk} \mathcal{A}' = r$. Let $B$ be a base of $\mathcal{A}'$ whence also a base of $\mathcal{A}$. Since $B$ is disjoint with $\mathcal{A}_X$ no $b \in B$ depends on $C$. □

**Corollary 6.1.** *Suppose for all* $X \in L(\mathcal{A}) \setminus \{0\}$ *we have*

$$|\mathcal{A}(X)| < \frac{n}{2}. \tag{1}$$

*Then* $\mathcal{A}$ *is large.*

Clearly it suffices to check the inequality (1) for $X$ of rank $r - 1$ only.

**Example**
The arrangements of the following classes of are clearly large.
(1) Generic arrangements with $|\mathcal{A}| \geq 2r - 1$.
(2) Every arrangement containing a large subarrangement of full rank.

## 7. Groups generated by reflections

**Definition 7.1.** *Let* $V$ *be a complex linear space of dimension* $r$. *A (complex) reflection is a finite order invertible linear transformation* $\tau : V \to V$ *whose fixed point set is a hyperplane (denoted* $H_\tau$). *A finite subgroup of* $GL(V)$ *is a reflection group if it is generated by reflections.*

For a reflection group $W$ the set $\mathcal{A}_W = \{H_\tau\}$ is called the *reflection arrangement* of $W$.

A reflection group $W$ is *irreducible* if its tautological representation to $GL(V)$ is irreducible. Then the rank of $W$ is $r$.

**Theorem 7.1.** *(see [3, 1, 5, 9]).* *Let* $M_W = V \setminus \bigcup_{H \in \mathcal{A}_W} H$ *for an arbitrary reflection group* $W$. *Then* $M_W$ *is a* $K[\pi, 1]$.

**Example**

For $\ell > 1$ every hyperplane $H_{ij} \subset \mathbb{R}^\ell$ of the Braid arrangement is the collection of fixed points of a real reflection permuting $x_i$ and $x_j$. Thus the (complexified) Braid arrangement is the reflection arangement for the permutation group $W = \Sigma_\ell$. Here $\pi_1(M_W)$ is the pure Braid group on $\ell$ strings, that is the *pure Artin group of type* $A_{\ell-1}$.

Similarly, for any (complexified) finite Coxeter group $W$ the group $\pi_1(M_W)$ is the pure Artin group of the respective type. Because of that $\pi_1(M_W)$ for an arbitrary finite complex reflection group $W$ is called the *pure Artin type group for* $W$ (or *the generalized pure Braid group associated to* $W$).

8. CALCULATION OF $TC_s(M_W)$

Here is the main theorem of the paper.

**Theorem 8.1.** *For every irreducible reflection group $W$ of rank $r$ and every $s > 1$ the arrangement $\mathcal{A}_W$ is well-balanced whence $TC_s(M_W) = sr - 1$.*

*Proof.* Our proof consists of four parts.

1. If $W$ has rank equal to 2 then the result is immediate since $r = 2$.

2. Infinite series.

For the infinite series well-balanced pairs can be exhibited explicitly. For that we give hyperplanes by their defining linear forms and $\mathcal{A}$ the product of them.

(a) **Full monomial types** $G(m, 1, r)$ (if $m = 1$ it is $\mathbf{A_r}$; if $m = 2$ it is $\mathbf{B_r}$) : $Q = \prod_{i=1}^{r} x_i \prod_{1 \leq i < j \leq r}(x_i^m - x_j^m)$. Put $B = \{x_1, \ldots, x_r\}$ and $C = \{x_1 - x_2, \ldots, x_1 - x_r\}$.

(b) **Special monomial types** $G(m, m, r)$ $m \geq 2$: $Q = \prod_{1 \leq i < j \leq r}(x_i^m - x_j^m)$ (if $m = 2$ it is $\mathbf{D_r}$). Put $B = \{x_1 - \zeta x_2, \ldots, x_1 - \zeta x_r, x_2 - \zeta x_3\}$ and $C = \{x_1 - x_2, \ldots, x_1 - x_r\}$ where $\zeta$ is a primitive root of 1 of order $m$.

In (a), the result is clear. In (b), $B$ is independent since it generates the basis $\{x_1, \ldots, x_r\}$ of $V^*$. Besides $C$ lies in the kernel of the index (the linear map $ind : V^* \to \mathbb{C}$, $ind(x_i) = 1$) while no $b \in B$ does.

3. The exceptional groups different from Coxeter types $\mathbf{E}_m$

In this case, we check case-by-case that $L(\mathcal{A}_W)$ is well-balanced using Tables C.1-C.23 from the book [10].

We use Corollary 6.28 from this book stating that $\mathcal{A}_X$ is the reflection arrangement for a reflection subgroup $W_X$ of $W$. The numbers $n_X = |\mathcal{A}_X|$ can be found from Table B.1 as the sums of covariants for $W_X$.

The table below is organized as follows. The first row consists of the Shephard-Todd classification numbers (23-34) of exceptional groups of ranks greater than 2 (no types $\mathbf{E}_m$). The second row consists of the cardinalities $n$ of the respective arrangements. The third row consists of the maximal cardinalities of $\mathcal{A}_X$. It suffices to check inequality (1): $|\mathcal{A}_X| < \frac{n}{2}$.

| 23 | 24 | 25 | 26 | 27 | 28 | 29 | 30 | 31 | 32 | 33 | 34 |
|----|----|----|----|----|----|----|----|----|----|----|----|
| 15 | 21 | 12 | 21 | 45 | 24 | 40 | 60 | 60 | 40 | 45 | 126 |
| 5 | 4 | 4 | 5 | 5 | 9 | 12 | 15 | 15 | 12 | 12 | 45 |

4. Types $\mathbf{E}_m$.

For the types $\mathbf{E}_m$, the inequality (1) does not hold but it is easy to check that $L(\mathcal{A})$ is well-balanced by definition. The needed information is in the table below.

| $E_6$ | $E_7$ | $E_8$ |
|-------|-------|-------|
| 36 | 63 | 120 |
| (20,15) | (36,21) | (63,42) |

The second row has the same meaning as in the previous table. The last row consists of pairs combining the maximal cardinality of $\mathcal{A}_Y$ with $\mathrm{rk}\, Y = r - 1$ and the cardinality of another $\mathcal{A}_X$ also with $\mathrm{rk}\, X = r - 1$. One needs to check that the sum in each pair is less than the entry of the second row. This shows that $L(\mathcal{A}_W)$ is well-balanced and finishes the proof.

□

**Remark.** Generic arrangements with $n < 2r - 1$ are not-large The general formula for generic arrangements is

$$\text{TC}_s(M) = \min\{sr - 1, (s-1)n\}.$$

For instance, if $r = 3$, $n = 4$, $s = 2$ then $\text{TC}_2(M) = 4 < 2r - 1$.

This result has been generalized and will be published in another paper.

**Conjecture**

For every complex hyperplane arrangement the topological complexity of its complement equals (for every $s$) the cohomological lower bound.

*E-mail address*: `yuz@math.uoregon.edu`

University of Oregon, Eugene, OR 97403